\documentclass[10pt,twoside,leqno]{article}

\usepackage{amsfonts}
\usepackage{amsmath}
\usepackage{amssymb}
\usepackage{amsopn}
\usepackage{amsthm}
\usepackage{amstext}

\renewcommand{\thefootnote}{}                                            
\newcounter{supersection}[section]                                       
\newtheorem{thm}[supersection]{Theorem}                                  
\newtheorem{lm}[supersection]{Lemma}                                     
\newtheorem{re}[supersection]{Remark}                                    
\newtheorem{ex}[supersection]{Example}                                   
\include{epsf}                                                           
\def\draw#1#2#3{\removelastskip                                          
\vskip 0cm plus #1cm                                                     
\goodbreak                                                               
\vskip 0cm plus -#1cm \vskip 0.3cm                                       
\vspace*{  0.0 mm} 
\hspace*{-10.0 mm} 
\dimen0=#2 true bp \ifdim\dimen0>\hsize\dimen0=\hsize\fi                 
\epsfxsize=\dimen0 \noindent {\epsfbox{#3.eps}} \count255=\epsfurx       
\advance\count255 by-\epsfllx \message{[#3.eps bp width =                
\number\count255]}                                                       
\goodbreak}                                                              
\def\bibname{\mbox{\normalsize \bf REFERENCES}}                          
\def\thebibliography#1{\paragraph*{\uppercase{\bibname}}\list            
{[\arabic{enumi}]}{\settowidth\labelwidth{[#1]}\leftmargin\labelwidth    
\advance\leftmargin\labelsep\usecounter{enumi}}                          
\def\newblock{\hskip .11em plus .33em minus .07em}                       
\sloppy\clubpenalty4000\widowpenalty4000                                 
\sfcode`\.=1000\relax}                                                   
\def\stop{\mbox{\footnotesize {\vrule width 6pt height 6pt}}}            

\begin{document}

\renewcommand{\thefootnote}{\fnsymbol{footnote}}

\setcounter{page}{1}

\pagestyle{myheadings}
\markboth{\rm  B. Male\v sevi\' c, M. Obradovi\' c}{\rm An application of Groebner bases to planarity of intersection of surfaces}

\date{ }

{ \footnotesize \noindent 

\noindent {Available at:} {\tt http://arxiv.org/} }

\bigskip

\leftline{} \hrule
\bigskip
\bigskip

\font\bigbf=cmbx12

\thispagestyle{empty} 

\font\bigbf=cmbx12

\centerline{\bigbf AN APPLICATION OF GROEBNER BASES TO}

\smallskip

\centerline{\bigbf PLANARITY OF INTERSECTION OF SURFACES}

\bigskip
\centerline{\bf Branko Male\v sevi\' c\footnote{First author supported in part by the project MNTRS, Grant No. ON144020.},
\bf Marija Obradovi\' c}

\footnote[0]{2000 {\it Mathematics Subject Classifications}.
14N05, 51N35.} \footnote[0]{{\it Key words and Phrases}.
Groebner bases, egg curve based conoid, planar intersection.}

\footnote[0]{
}

\footnote[0]{}


\begin{abstract}
In this paper we use Groebner bases theory in order to determine planarity of intersections
of two algebraic surfaces in $\mbox{\bf R}^3$. We specially considered plane sections of certain
type of conoid which has a cubic egg curve as one of the directrices. The paper investigates
a possibility of conic plane sections of this type of conoid.
\end{abstract}

\section{\mbox{\normalsize \bf 1. THE BASIC CONCEPT OF GROEBNER BASES}}

\medskip
\noindent
Many problems in mathematics can be treated as a problem of solvability of a system of polynomial equations:
\begin{equation}
\label{EQ_01}
\begin{array}{c}
f_{1}(x_1, x_2, \ldots, x_n) = 0 \\
f_{2}(x_1, x_2, \ldots, x_n) = 0 \\
\vdots                           \\
f_{m}(x_1, x_2, \ldots, x_n) = 0
\end{array}
\end{equation}
where $f_{1}, f_{2}, \ldots, f_{m}$ are polynomials in $n$ variables over set of complex numbers \mbox{\bf C}.
The main matter in the theory of systems of polynomial equations are Groebner bases.
Let us define this concept through following considerations. The polynomial of several variables:

\vspace*{-1.5 mm}

\begin{equation}
\label{EQ_02}
P(x_1,\ldots,x_n)
=
\displaystyle\sum\limits_{\alpha \in A_{0}}{c_{\alpha} \, x_{1}^{\alpha_1} \cdots x_{n}^{\alpha_n}}
\end{equation}
is considered as a sum over $n$-tuples  $\alpha\!=\!(\alpha_{1}, \ldots, \alpha_{n})$ from a finite set
\mbox{$\!A_{0}\!\subset\!\mbox{\bf N}_{0}^{n}\!$}~\mbox{$(c_{\alpha}\!\neq\!0)$}. The term $x_{1}^{\alpha_1} \cdots x_{n}^{\alpha_n}$
is called {\em the monomial}, while the summand $c_{\alpha} \, x_{1}^{\alpha_1} \cdots x_{n}^{\alpha_n}$ is called
{\em the term of the polynomial}.

\break

\noindent
A {\em monomial order} $\succ$ is a total order on the set of all monomials, satisfying two additional properties:
the ordering respects multiplication (if $u \succ v$ and $w$ is any other monomial, then $uw \succ vw$) and
the ordering is a well ordering. {\em The lexicographic order} over a set of monomials~$\succ_{lex}$ is defined by
{\em order of variables}
$
x_{1} \succ_{lex} \ldots \succ_{lex} x_{n}
$
so that:
$
x_{1}^{\alpha_{1}} \cdots \, x_{n}^{\alpha_{n}}
\succ_{lex}
x_{1}^{\beta_{1}} \cdots \, x_{n}^{\beta_{n}}
\; \mbox{iff} \;
\alpha_{1}\!=\!\beta_{1}, \ldots, \alpha_{k-1}\!=\!\beta_{k-1}, \alpha_{k} \!>\! \beta_{k}
$
for some $k$. Other types of monomial ordering are presented in \cite{[Froberg_1997]}.
Let~us notice that polynomial (\ref{EQ_02}) can be {\em arranged} in the given monomial order, so  that the first term
in the sum $f$ is {\em the leading term} and denoted by $LT(f)$. Then leading term is shown as the product
of {\em the leading coefficient} $LC(f)$ and {\em the leading monomial} $LM(f)$, i.e. $LT(f) = LC(f) \cdot LM(f)$.

\medskip
\noindent
In the polynomial ring $\mbox{\bf C}[x_{1}, \ldots, x_{n}]$ we determine a {\em division algorithm}, which follows.
Let a $m$--tuple of polynomials $F = (f_{1}, \ldots, f_{m})$ be given, arranged into monomial order $\succ$ over the set of monomials.
Then we quote an division algorithm with the following pseudo-code \cite{[IVA_2007]}:
\begin{equation}
\label{EQ_03}
\mbox{\footnotesize \sf
\begin{tabular}{llll}
Input:  & \multicolumn{3}{l}{$f, f_{1}, \ldots, f_{m}$}                                                                                \\[+0.1 ex]
Output: & \multicolumn{3}{l}{$r, a_{1}, \ldots, a_{m}$}                                                                                \\[+0.9 ex]
        & \multicolumn{3}{l}{$r:=0, a_{1}:=0, \ldots, a_{m}:=0 $}                                                                      \\[+0.1 ex]
        & \multicolumn{3}{l}{$p:=f$}                                                                                                   \\[+0.1 ex]
        & \multicolumn{3}{l}{{\sc while} $p \neq 0$ {\sc do}}                                                                          \\[+0.1 ex]
        & \multicolumn{3}{l}{\quad \qquad $i := 1$}                                                                                    \\[+0.1 ex]
        & \multicolumn{3}{l}{\quad \qquad {\rm divisionoccured} := {\rm false}}                                                        \\[+0.1 ex]
        & \multicolumn{3}{l}{\quad \qquad {\sc while} $\; i \leq m\; $ {\sc and} {\rm divisionoccured} = {\rm false} {\sc do}}         \\[+0.1 ex]
        & \hspace*{9.0 mm} & \multicolumn{2}{l}{\quad \quad \enskip {\sc if} {\it LT}($f_{i}$) {\rm divides} {\it LT}($p$) {\sc then}} \\[+0.1 ex]
        & \hspace*{9.0 mm} & \quad \qquad \quad $a_{i} := a_{i} + \mbox{\it LT}(p)/\mbox{\it LT}(f_{i})$                               \\[+0.1 ex]
        & \hspace*{9.0 mm} & \quad \qquad \quad $p := p - \mbox{\it LT}(p)/\mbox{\it LT}(f_{i}) \cdot f_{i}$                           \\[+0.1 ex]
        & \hspace*{9.0 mm} & \quad \qquad \quad {\rm divisionoccured} := {\rm true}                                                    \\[+0.1 ex]
        & \hspace*{9.0 mm} & \multicolumn{2}{l}{\quad \quad \enskip \sc else}                                                          \\[+0.1 ex]
        & \hspace*{9.0 mm} & \quad \qquad \quad $i := i + 1$                                                                           \\[+0.1 ex]
        & \multicolumn{3}{l}{\quad \qquad \enskip {\sc if} {\rm divisionoccured} = {\rm false} {\sc then}}                             \\[+0.1 ex]
        & \hspace*{9.0 mm} & \multicolumn{2}{l}{\quad $r := r + \mbox{\it LT}(p)$}                                                     \\[+0.1 ex]
        & \hspace*{9.0 mm} & \multicolumn{2}{l}{\quad $p := p - \mbox{\it LT}(p)$}                                                     \\[+0.1 ex]
        & (\mbox{\sc stop}) &
\end{tabular}}
\end{equation}

\noindent
Using the previous algorithm, every polynomial $f \in \mbox{\bf C}[x_{1}, \ldots, x_{n}]$ can be presented as follows:

\vspace*{-3.0 mm}

\begin{equation}
\label{EQ_04}
f = a_{1} f_{1} + \ldots + a_{m} f_{m} + r
\end{equation}

\noindent
for some $a_{i}, r \in \mbox{\bf C}[x_{1}, \ldots, x_{n}]$ $(i = 1, \ldots, m)$, while either $r = 0$ or $r$ is $\mbox{\bf C}$-linear
combination of monomials none of which is divisible by any $LT(f_{1}), \ldots, LT(f_{m})$. The polynomial $r=REM(f,F)$ is called
{\em the remainder} of division of $f$ by \mbox{$m$--tuple} of polynomials $F$ arranged into monomial order $\succ$. The main characteristic
of the remainder $r$ is that it is not uniquely determined in comparison to order of division by polynomials from  $m$--tuple $F$.
Let us remark:
\begin{equation}
\label{EQ_05}
LT(f) \succ LT(r).
\end{equation}
For a $m$--tuple of polynomials  $F = (f_1, \ldots, f_m)$ {\em the ideal} is  determined by:
\begin{equation}
\label{EQ_06}
I
=
\left\{
\displaystyle\sum\limits_{i=1}^{m}{h_{i} f_{i} \; \mbox{\Large $|$} \; h_{1}, \ldots, h_{m} \in \mbox{C}[x_{1}, \ldots, x_{n}]}
\right\}.
\end{equation}

\noindent
It is denoted by $\langle f_{1}, \ldots, f_{m} \rangle$. Then polynomials $f_{1}, \ldots, f_{m}$ are called a {\em generators} for
the ideal $I = \langle f_{1}, \ldots, f_{m} \rangle$ in the polynomial ring $\mbox{\bf C}[x_{1}, \ldots, x_{n}]$. Let us emphasize
that Hilbert Basis Theorem claims that every ideal has a finite generating set \cite{[IVA_2007]}. Next, for the ideal
$I = \langle f_{1}, \ldots, f_{m} \rangle$ we define {\em the ideal of leading terms}:
\begin{equation}
\label{EQ_07}
\langle
LT(I)
\rangle
=
\langle
\{
LT(f) \, | \, f \!\in\! I
\}
\rangle .
\end{equation}
In any ideal $I = \langle f_{1}, \ldots, f_{m} \rangle$ it is true:
\begin{equation}
\label{EQ_08}
\langle
LT(f_{1}), \ldots, LT(f_{m})
\rangle
\subseteq
\langle LT(I)
\rangle
\end{equation}
Generating set $G = \{g_1, \ldots, g_s \}$ of the ideal $I = \langle f_1, \ldots, f_m \rangle$ is called a {\em Groebner basis} if

\vspace*{-3.0 mm}

\begin{equation}
\label{EQ_09}
\langle
LT(g_{1}), \ldots, LT(g_{s})
\rangle
=
\langle LT(I)
\rangle
\end{equation}
Let us remark that for each permutation of $s$--tuple $G = \langle g_{1}, \ldots, g_{s} \rangle$ the remainder $r = REM(f, G)$
is uniquely determined \cite{[IVA_2007]}.

\medskip
\noindent
{\sc B. Buchberger} in his dissertation has shown that for each ideal there exists a Groebner basis by following algorithm
\cite{[IVA_2007]}:
\begin{equation}
\label{EQ_10}
\mbox{\footnotesize \sf
\begin{tabular}{llll}
Input:  & \multicolumn{3}{l}{$F = \{f_{1}, \ldots, f_{m}\}$--{\rm Set of generators of} $I$}                               \\[+0.1 ex]
Output: & \multicolumn{3}{l}{$G = \{g_{1}, \ldots, g_{s}\}$\,--$\;${\rm Groebner basis for} $I$ $(G \supseteq F)$}         \\[+0.4 ex]
        & \multicolumn{3}{l}{$G:=F$}                                                                                       \\[+0.1 ex]
        & \multicolumn{3}{l}{\sc repeat}                                                                                   \\[+0.1 ex]
        & & \hspace*{-5.0 mm}   {$G\,':=G$}                                                                                \\[+0.4 ex]
        & & \hspace*{-5.0 mm}   {\sc for} {\rm each pair} $\{p,q\}$, $p \neq q$, {\rm in} $G\,'$ {\sc do}                  \\[+0.1 ex]
        & & \hspace*{2.0 mm}    $S:=S(p,q)$                                                                                \\[+0.1 ex]
        & & \hspace*{2.0 mm}    $h:=REM(S,G\,')$                                                                           \\[+0.1 ex]
        & & \hspace*{2.0 mm}    {\sc if} $h \neq 0$ {\sc then} $G:=G \cup \{h\}$                                           \\[+0.1 ex]
        & & \hspace*{-13.5 mm}  {\sc until} $G = G\,'$                                                                     \\[+0.1 ex]
        & (\mbox{\footnotesize \sc stop}) &
\end{tabular}}
\end{equation}
In the previous algorithm we used {\em the $S$-polynomial} of $f$ and $g$ (in relation to the fixed monomial order
$\succ$) which is defined by:
\begin{equation}
\label{EQ_11}
S
=
S(f,g)
=
\displaystyle\frac{\mbox{\small $LCM$}(\mbox{\small $LM$}(f),\mbox{\small $LM$}(g))}{\mbox{\small $LT$}(f)} \cdot f
-
\displaystyle\frac{\mbox{\small $LCM$}(\mbox{\small $LM$}(f),\mbox{\small $LM$}(g))}{\mbox{\small $LT$}(g)} \cdot g,
\end{equation}
where $LCM$ is the least common multiple. For different sets of generators of the ideal~$I$ Groebner basis,
produced by {\em the Buchberger's algorithm}, is not unique. {\em Groebner basis is minimal} iff $LC(g) = 1$
for all $g \in G$, whereat $\{LT(g) \, | \, g \!\in\! G\}$ is a minimal basis of the monomial ideal $\langle \{LT(f) \, | \, f \!\in\! I\} \rangle$
(see Exercise 6/\S7/Ch.$\,$2~\cite{[IVA_2007]}). Minimal Groebner basis is not unique as well. {\em Reduced Groebner basis}
is minimal Groebner basis for the ideal $I$ such that for all $p \in G$ there is no monomial of $p$ which lies in
$\langle \{LT(g) \, | \, g \!\in\! G \backslash \{p\} \,\} \rangle $ \cite{[Froberg_1997]}, \cite{[IVA_2007]}.
The following statement is true \cite{[IVA_2007]}:

\medskip
\noindent
\begin{thm}\label{Theorem_1_1}
Every polynomial ideal has a unique reduced Groebner basis, related to a fixed monomial order.
\end{thm}

\medskip
\noindent
Groebner bases in applications are usually considered in reduced form, as it is emphasized in
\cite{[IVA_2007]} and \cite{[Roanes-Lozano_2004]}. In the next part of this paper we consider an
application of the Groebner bases in the real three-dimensional geometry.

\section{\mbox{\normalsize \bf 2. AN APPLICATION OF GROEBNER BASES TO PLANARITY}}

\vspace*{-2.5 mm}

\qquad\enskip\enskip{\mbox{\normalsize \bf OF INTERSECTION OF SURFACES}}

\bigskip
\noindent
Let us consider the consistent system of two real non-linear polynomial equations with three variables:

\vspace*{-3.0 mm}

\begin{equation}
\label{EQ_12}
f_{1}(x,y,z)=0 \quad \wedge \quad f_{2}(x,y,z) = 0.
\end{equation}
For the system $(\ref{EQ_12})$ we define that {\em the system has planar solution} if there exists a linear-polynomial:

\vspace*{-3.0 mm}

\begin{equation}
\label{EQ_12'}
g=g(x,y,z)=Ax + By + Cz + D,
\end{equation}
for some real constants $A,B,C,D$, such that every solution of the system (\ref{EQ_12})
is also a solution of the linear equation:
\begin{equation}
\label{EQ_13}
g(x,y,z) = 0.
\end{equation}
Thus we indicate that {\em the system has planar intersection}. Let a monomial order be determined, then the following statement is true:

\begin{thm}\label{Theorem_2_1}
If a Groebner basis of the system $(\ref{EQ_12})$ contains the linear polynomial $(\ref{EQ_12'})$,
then the solution of the system $(\ref{EQ_12})$ is planar.
\end{thm}

\medskip
\noindent
The following statement is formulated for lexicographic order \mbox{$x \succ_{lex} y \succ_{lex} z$}:

\begin{thm}\label{Theorem_2_2}
Let the system $(\ref{EQ_12})$ have the planar intersection by$:$
\begin{equation}
\label{EQ_14}
Ax + By + Cz + D = 0,
\end{equation}
for some $A, B, C, D \in \mbox{\bf R}$ and $A \neq 0$. If for the ideal $I=\langle f_{1}, f_{2} \rangle$
the following is true$:$

\vspace*{-3.0 mm}

\begin{equation}
\label{EQ_14'}
y \not\in \langle LT(I) \rangle
\quad \wedge \quad
z \not\in \langle LT(I) \rangle,
\end{equation}
then the linear polynomial$:$
\begin{equation}
\label{EQ_15}
\hat{g}
=
\hat{g}(x,y,z)
=
x + (B/A) \, y + (C/A) \, z + (D/A)
\end{equation}
is an element of the reduced Groebner basis related to lexicographic order \mbox{$\succ_{lex}$}.
\end{thm}

\medskip
\noindent
{\bf Proof.} The system $(\ref{EQ_12})$ is equivalent with the following system:
$$
f_{1}(x,y,z)=0 \quad \wedge \quad f_{2}(x,y,z) = 0 \quad \wedge \quad \hat{g}(x,y,z) = 0.
\leqno (\ref{EQ_12}')
$$
The reduced Groebner bases for systems $(\ref{EQ_12})$ and $(\ref{EQ_12}')$ are equal. Let us prove that $\hat{g} \in I$ is element
of reduced Groebner basis for system $(\ref{EQ_12}')$. Using Buchberger algorithm, we can assume that polynomial $\hat{g}$ is element
of the Groebner basis $G_{1}$. Let~us consider minimal Groebner basis $G_{min}$ formed from $G_{1}$ such that
$\hat{g} \in G_{min}$. Next, let~us prove that polynomial $\hat{g}$ is reduced related to $G_{min}$
(see~Proposition~6./\S7/Ch.$\,$2 in \cite{[IVA_2007]} or Lect.$\,$14 in \cite{[Karp_2006]}).
By condition (\ref{EQ_14'}) it follows that for a linear polynomial $\hat{g} = x + (B/A)y + (C/A)z + (D/A)$
it is true that $x, y, z \not \in \langle LT {\big (} G_{min} \backslash \{\hat{g}\} {\big )} \rangle$.
If~polynomial $\hat{g}$ is reduced in one minimal Groebner basis, then $\hat{g}$ is reduced related to every minimal
Groebner basis (see proof of the Proposition~6./\S7/Ch.$\,$2 in \cite{[IVA_2007]}) and therefore polynomial $\hat{g}$
is an element of the unique reduced Groebner basis.$\;$\stop

\smallskip\noindent
Let us emphasize that computer algebra systems compute reduced Groebner basis. In this paper we use computer algebra system Maple
and Maple Package "{\sf Groebner}" which is initiated by:
$$
\mbox{\sf with(Groebner);}
$$
Let $\mbox{\sf F} := [\, \mbox{\sf f}_{1}, \mbox{\sf f}_{2} \,]$ be the list of two polynomials
of three variables $\mbox{\sf x}, \mbox{\sf y}, \mbox{\sf z}$. In Maple lexicographic order
$\mbox{\sf x} \succ_{lex} \mbox{\sf y} \succ_{lex} \mbox{\sf z}$ we denote by $\mbox{\sf plex}(\mbox{\sf x}, \mbox{\sf y}, \mbox{\sf z})$.
Then by command
$$
\mbox{\sf Basis}(\mbox{\sf F},\mbox{\sf plex}(\mbox{\sf x}, \mbox{\sf y}, \mbox{\sf z}));
$$
Maple computes the reduced Groebner basis. If some variables are omitted from the list of variables, then the algorithm considers these
variables as non-zero constants (see Appendix C/\S2\,({\sf Maple}) in \cite{[IVA_2007]} and \cite{[Heck_1997]}, \cite{[Heck_2003]}).

\begin{ex}\label{Example_2_3.}
Let the system of two polynomials be given$:$
\begin{equation}
\label{EQ_16}
\mbox{\sf F}
:=
{\bigg [}
\mbox{\sf z} - \displaystyle\frac{\mbox{\sf x}^2}{\mbox{\sf a}^2} - \displaystyle\frac{\mbox{\sf y}^2}{\mbox{\sf b}^2},
\displaystyle\frac{\mbox{\sf x}^2}{\mbox{\sf a}^2} + \displaystyle\frac{\mbox{\sf y}^2}{\mbox{\sf b}^2}
- \displaystyle\frac{\mbox{\sf x}}{\mbox{\sf a}} - \displaystyle\frac{\mbox{\sf y}}{\mbox{\sf b}}
{\bigg ]};
\end{equation}
where $\mbox{\sf x},\mbox{\sf y},\mbox{\sf z}$ are variables and $\mbox{\sf a}, \mbox{\sf b}$ are non-zero constants.
We compute reduced Groebner basis using Maple by command
$
\mbox{\sf Basis}(\mbox{\sf F}, \mbox{\sf plex}(\mbox{\sf x},\mbox{\sf y},\mbox{\sf z}));
$
and the result is the following list of the $($arranged$)$ polynomials$:$
\begin{equation}
\label{EQ_17}
{\big [}
\mbox{\sf bx} + \mbox{\sf ay} - \mbox{\sf abz}, 2\mbox{\sf y}^2-2\mbox{\sf byz} - \mbox{\sf b}^2\mbox{\sf z}^2 - \mbox{\sf b}^2\mbox{\sf z}
{\big ]}
\end{equation}
Let us emphasize that the answer is the reduced Groebner basis, except for clearing denominators $($so leading coefficients
of the Groebner basis are polynomials~in $\mbox{\sf a}$~and~$\mbox{\sf b}$$)$ {\rm \cite{[IVA_2007]}}.
This proves that curve of intersection lies in the following~plane$:$
\begin{equation}
\label{EQ_18}
\mbox{\sf bx} + \mbox{\sf ay} - \mbox{\sf abz} = \mbox{\sf 0}.
\end{equation}
Let us remark that $\langle LT(I) \rangle = \langle \mbox{\sf bx}, 2\mbox{\sf y}^2 \rangle$ and the condition $(\ref{EQ_14'})$ is fulfilled.
\end{ex}

\begin{ex}\label{Example_2_4}
Let the system of two polynomials be given$:$
\begin{equation}
\label{EQ_16'}
\mbox{\sf F}
:=
{\big [}
\mbox{\sf x} + \mbox{\sf yz} + \mbox{\sf y} - \mbox{\sf z}^4 - 4, \mbox{\sf y} - \mbox{\sf z}^3 - 1
{\big ]};
\end{equation}
where $\mbox{\sf x},\mbox{\sf y},\mbox{\sf z}$ are variables. We compute reduced Groebner basis using Maple by command
$
\mbox{\sf Basis}(\mbox{\sf F}, \mbox{\sf plex}(\mbox{\sf x},\mbox{\sf y},\mbox{\sf z}));
$
and the result is the following list of the $($arranged$)$ polynomials$:$

\vspace*{-3.0 mm}

\begin{equation}
\label{EQ_17'}
{\big [}
\mbox{\sf x} + \mbox{\sf z}^3 + \mbox{\sf z} - 3, \mbox{\sf y} - \mbox{\sf z}^3 - 1
{\big ]}
\end{equation}
Note that$:$

\vspace*{-3.0 mm}

\begin{equation}
\label{19'}
1 \!\cdot\! (\mbox{\sf x} + \mbox{\sf yz} + \mbox{\sf y} - \mbox{\sf z}^4 - 4) - \mbox{\sf z} \!\cdot\! (\mbox{\sf y} - \mbox{\sf z}^3 - 1)
=
\mbox{\sf x} + \mbox{\sf y} + \mbox{\sf z} - 4
\end{equation}
and therefore the intersection is planar. Let us emphasize that $\langle LT(I) \rangle = \langle \mbox{\sf x}, \mbox{\sf y} \rangle$
and condition $(\ref{EQ_14'})$ is not fulfilled.
\end{ex}

\begin{re}\label{Remark_2_5.}
It is possible to formulate and prove the previous theorems in more extended sense than in the given three-dimensional formulations
and lexicographic order.
\end{re}

\section{\mbox{\normalsize \bf 3. AN  APPLICATION OF THE GROEBNER BASES TO PLANE}}

\vspace*{-2.5 mm}

\qquad\enskip\enskip{\mbox{\bf SECTIONS OF ONE TYPE OF CONOID WITH BASIC CUBIC}}

\vspace*{1.5 mm}

\quad\enskip{\mbox{\bf EGG CURVE}}

\medskip
\noindent
In this part, there will be considered an application of the Groebner bases theory on a surface which is obtained
by motion of the system of generatrices along three directrices.

\smallskip
\noindent
Let us accept cubic egg curve $g$, the right cubic hyperbolic parabola of type $A$ \cite{[Dovnikovic_1]} which lies
in the plane $x-y$, for the initial curve taken as a plane section of a conoid, and also for its directrix. Starting
from the definition of conoid \cite{[Rovenskii_4]} (p.$\,$277) as a surface which has two directrices in finiteness:
plane curve $d_1$, and straight line $d_2$,  while the third directrix, straight line $d_3$, lies in infinity,
we consider a conoid with \cite{[Obradovic_Petrovic_9]}:


\smallskip
\noindent
$(i)$ Cubic curve $(g)$ in the form $g=d_1$: $b^2 x^2 + a^2 y^2 + 2 d x y^2 + d^2 y^2 - a^2 b^2 = 0 $ as directrix $d_1$
($a\!>\!b\!>\!0,\;a-b \geq d > 0$).

\smallskip
\noindent
$(ii)$ Straight line parallel to the axis $y$, in the plane $y-z$ on elevation $z=h$ as directrix $d_2$ $(h\!>\!0)$.

\smallskip
\noindent
$(iii)$ Plane $x-z$ as directrix plane of the conoid, and its infinitely distant straight line as directrix $d_3$,
in order to obtain a right conoid.


\break

\vspace*{115.0 mm} 

\smallskip

\centerline{Figure 1.}

\bigskip

\noindent
Each generatrix of the surface will intersect all three directrices \cite{[Sbutega_Zivanovic_5]}.
The surface obtained in this manner will be quintic surface \cite{[Obradovic_Petrovic_9]},
degenerated into a conoid of the fourth order, as shown at Figure 1, with the equation \cite{[Obradovic_Petrovic_9]}:
\begin{equation}
\label{EQ_19}
(a^2 y^2 + d^2 y^2 - a^2 b^2)(z-h)^2 - 2dhxy^2(z - h) + b^2 h^2 x^2 = 0
\end{equation}
and a plane:

\vspace*{-3.0 mm}

\begin{equation}
\label{EQ_20}
z=h.
\end{equation}
Thus, we give an explanation of origination and the order of the surface, from the aspect of Projective Geometry.

\smallskip
\noindent {\boldmath $1.$}
As every plane $F_1$ parallel to the plane $x-z$, which will imply generatrices of the surface, intersects the plane $x-y$
of the curve $g$ by a straight line parallel to the $x$ axes \cite{[Gray_2]}, and each  straight line must intersect
the cubic curve by three points $X_1$, $X_2$ and $X_3$, we conclude that in each point of  the directrix $d_2$ there
intersect three generatrices: $i_1$, $i_2$, and $i_3$, of the surface. Therefore, the directrix $d_2$ will be
triple line of the surface \cite{[Gorjanc_3]}.

\smallskip
\noindent {\boldmath $2.$}
Because of the symmetry of the curve $g$ in relation to the $x$ axes, pair of generatrices of the conoid will intersect
in the same point of the infinite directrix $d_3$, so it  will be double line of the surface.

\smallskip
\noindent {\boldmath $3.$}
Cubic curve $g=d_1$ has two asymptotes $as_1$ and $as_2$, whereat the asymptote $as_1$ is parallel to the $y$ axes,
and $as_2$ is parabola \cite{[Wilson_7]} with $x$ axes. Therefore, the cubic curve $g$ has two infinite points,
$X_{\infty}$  and $Y_{\infty}$, of the axes $x$ and $y$.

\smallskip
\noindent {\boldmath $4.$}
$X_{\infty}$ represents the infinite vertex of the parabola $as_2$. Every generatrix $i_3$ of the surface passes through
this point forming the plane $H_1$, jointly with the directrix $d_2$. $Y_{\infty}$  is a point of tangency of curve
$g = d_1$ and asymptote $as_1$. At the same time, it is a triple point of intersection of the directrix $d_2$ and the plane $x-y$.

\smallskip
\noindent {\boldmath $5.$}
Points $X_{\infty}$ and $Y_{\infty}$ belong to infinite straight line $q_{\infty}$ of the plane $x-y$, which passes two times through
those points. Straight line $q_{\infty}$ represents two overlapped generatrices $i_{2}$ and $i_{3}$ which pass through the infinite point
$Y_{\infty}$ and are tangents to the vertex $X_{\infty}$ of the parabolic asymptote $as_2$ (and also the curve $g$) at two immeasurable
close points, $X_2$ and $X_3$, which are collinear to the point $Y_{\infty}$. Therefore, generatrices $i_2(X_{2\infty}Y_{\infty})$
and $i_3(X_{3\infty}Y_{\infty})$ will be identical, and so the straight line $q_{\infty}$ will be also the double line. Generatrix
$i_1$ will be separated, so in the infinite plane $\omega\infty$, there will lie five lines: double line $d_{3\infty}$,
double line $q_{\infty}$, and the line $i_{1\infty}$.

\smallskip
\noindent {\boldmath $6.$}
Plane $x-y$ intersects the conoidal surface by the initial cubic curve $g$ and by the infinite double straight line $q_{\infty}$,
which indicates that it is the case of the surface of fifth degree \cite{[Obradovic_Petrovic_9]}.

\smallskip
\noindent {\boldmath $7.$}
Considering that the generatrix $i_3$ is a straight line that belongs to the pencil of straight lines through the point $X_{\infty}$,
which constitutes, along with directrix $d_2$, the plane $H_1$, the surface will be degenerated to the conoid and the plane $H_1$.
By eliminating the plane, the doubleness of the line $q_{\infty}$ subdues, as well as the tripleness of the directrix $d_2$,
so the conoidal surface remains of the fourth order.

\smallskip
\noindent {\boldmath $8.$}
The horizontal plane $H_1$, at the elevation $z=h$ of the directrix $d_2$ is asymptotic plane of the surface. It is tangential plane
onto the infinite double straight line $q_{\infty}$. Besides, there will exist yet another asymptotic plane, a skew plane $R$,
which is determined by the directrix $d_{2}$ and the asymptote $as_1$ of the curve $g=d_1$, as two parallel lines, and is tangential
the conoid by the infinite straight line $r_{\infty} = i_1$.

\smallskip
\noindent {\boldmath $9.$}
The infinite plane of the space $\omega_{\infty}$, intersects the conoid by one double straight line $d_3$ and also by the remaining
straight lines $q_{\infty}$, and $r_{\infty}$, which proves that it is the surface of fourth order.

\medskip
\noindent
Mathematical explanation of the originating surface, initially of the fifth order is easily obtained
by simple multiplication of the equation of the fourth order (\ref{EQ_19}) by factor $z-h$. The equation derived
in this manner:
$
(a^2 y^2 + d^2 y^2 - a^2 b^2)(z-h)^3 - 2dhxy^2\mbox{$(z - h)^2$} + b^2 h^2 x^2 (z-h) = 0,
$
decomposes to an equation of the fourth order (\ref{EQ_19}) and an equation of a plane (\ref{EQ_20}), as it is given
in the commentary of the Projective Geometry.

\medskip
\noindent
Let us consider if it is possible for the plane section of the surface (\ref{EQ_19}) to be a non-degenerated conic
\cite{[Obradovic_Petrovic_9]}. Let us emphasize that this problem is equivalent to the following system:
\begin{equation}
\label{EQ_21}
\left\{
\begin{array}{c}
\!(\ref{EQ_21}/i)\;\;
(a^2 y^2 + d^2 y^2 - a^2 b^2)(z-h)^2 - 2dhxy^2(z - h) + b^2 h^2 x^2 = 0, \\[2.0 ex]
\!(\ref{EQ_21}/ii)\;\;
Ax + By + Cz + D = 0;
\end{array}
\right\}
\end{equation}
which has solution by non-degenerated conic, for some real constants $A,B,C,D$. Using technique of Groebner bases we will prove that is not possible.
The following statements are true:

\begin{lm}\label{Lemma_3_1}
{\boldmath $(i)$}
Plane $x\!=\!\alpha$, for $\alpha \in \mbox{\bf R} \backslash \{0\}$, has the intersection with the surface
$(\ref{EQ_21}/i)$ by the curve of the fourth order. Plane $x\!=\!0$ has the intersection with the surface $(\ref{EQ_21}/i)$
by $x\!=\!0$, $z\!=\!h\,{\big (}|y|\!\leq\!b \vee |y|\!\geq\!ab/d{\big )}$ or $x\!=\!0$, $y\!=\!\pm {ab}/{\sqrt{a^2+d^2}}$.

\smallskip\noindent
{\boldmath $(ii)$} Plane $y\!=\!\beta$, for $\beta \in \mbox{\bf R}$, has intersection with the surface $(\ref{EQ_21}/i)$  iff$:$
\begin{equation}
\label{EQ_22}
\hspace*{33.0 mm}
|\beta| \leq b
\quad\vee\quad
|\beta| \geq \displaystyle\frac{ab}{d}
\hspace*{33.0 mm}
{\Big (}\displaystyle\frac{ab}{d} > b{\Big )};
\end{equation}
then the intersection is presented as two straight lines $($which determine degenerated conic$)$.

\smallskip\noindent
{\boldmath $(iii)$} Plane $z\!=\!\gamma$, for $\gamma \in \mbox{\bf R} \backslash \{h\}$, has the intersection with the surface
$(\ref{EQ_21}/i)$ by the curve of the third order. Plane $z\!=\!h$ has the intersection with the surface $(\ref{EQ_21}/i)$
by $x\!=\!0$, $z\!=\!h\,{\big (}|y|\!\leq\!b \vee |y|\!\geq\!ab/d{\big )}$.
\end{lm}

\smallskip
\noindent
{\bf Proof.} {\boldmath $(i)$} By substitution $x\!=\!\alpha \in \mbox{\bf R}$ in $(\ref{EQ_21}/i)$ we obtained:
$$
\begin{array}{rl}
&
{\big (}y^2(a^2+d^2) - a^2b^2{\big )}(z-h)^2 - 2 h d \alpha y^2(z-h) + h^2 b^2 \alpha^2 = 0           \\[1.0 ex]
\Longleftrightarrow \!\!\! &
{\big (}z-h{\big )}{\big (}a^2b^2 - (a^2+d^2)y^2{\big )}
=
{\Big (}\!-dy^2 \pm \displaystyle\sqrt{{\big (}b^2-y^2{\big )}{\big (}a^2b^2 - d^2 y^2{\big )}} \, {\Big )} h \alpha \,.
\end{array}
$$
Based on $a,d > 0$ we can conclude that the intersection is one algebraic curve of the forth order if $\alpha \!\neq\! 0$.
If $\alpha\!=\!0 \,(\, \Longleftrightarrow x\!=\!0 \,)$, then from the previous equation follows
$x\!=\!0$, $z\!=\!h\,{\big (}|y|\!\leq\!b \vee |y|\!\geq\!ab/d{\big )}$
or $x\!=\!0$, $y=\pm {ab}/{\sqrt{a^2+d^2}}$.

\smallskip\noindent
{\boldmath $(ii)$} By substitution $y\!=\!\beta \in \mbox{\bf R}$ in $(\ref{EQ_21}/i)$ we obtained:
$$
\begin{array}{l}
b^{2}\,h^{2}\,x^{2}
+
(a^{2}\,\beta ^{2}
+
d^{2}\,\beta ^{2}
-
a^{2}\,b^{2})\,z^{2}
-
2\,d\,h\,\beta ^{2}\,x\,z
+                                                                                                 \\[1.5 ex]
2\,d\,h^{2}\,\beta^{2}\,x
+
(2\,a^{2}\,b^{2}\,h
-
2\,d^{2}\,h\,\beta^{2}
-
2\,a^{2}\,h\,\beta^{2})\,z
+                                                                                                 \\[1.5 ex]
(
a^{2}\,h^{2}\,\beta^{2}
+
d^{2}\,h^{2}\,\beta^{2}
-
a^{2}\,b^{2}\,h^{2}
)
=
0.
\end{array}
$$
If we solve the previous second order equation by $x$, from discriminant, we obtain the condition
(\ref{EQ_22}) for existence of the intersection via two straight lines:
$$
x
=
\pm
\displaystyle\frac{\,\beta^2d + \displaystyle\sqrt{{\big (}b^2-\beta^2{\big )}{\big (}a^2b^2-d^2\beta^2{\big )}}\,}{b^2h}
\,
{\big (}\!-\!z+h {\big )}.
$$
{\boldmath $(iii)$}  By substitution $z\!=\!\gamma \in \mbox{\bf R}$ in $(\ref{EQ_21}/i)$ we obtained:
$$
\begin{array}{rl}
&
b^{2}\,h^{2}\,x^{2}
+
2\,d\,h\,(h-\gamma)^2\,x \, y^{2}
+
(a^2+d^2)(h-\gamma) \, y^2
-
a^2 \, b^2 (h-\gamma)^2
=
0                                                                                                 \\[1.0 ex]
\Longleftrightarrow \!\!\! &
x
=
\pm
\displaystyle\frac{\displaystyle\sqrt{{\big (}b^2-y^2{\big )}{\big (}a^2b^2-d^2y^2{\big )}}
- y^2d}{b^2h}\,{\big (}h - \gamma{\big )}\,.
\end{array}
$$
Based on $d,h > 0$ we can conclude that the intersection is one algebraic curve of the third order iff
$\gamma \!\neq\! h$. If $\gamma\!=\!h \,(\, \Longleftrightarrow z\!=\!h \,) $, then from the previous equation
follows $z\!=\!h$, $x\!=\!0\,{\big (}|y|\!\leq\!b \vee |y|\!\geq\!ab/d{\big )}$.$\;$\stop
\begin{re}\label{Remark_3_2}
Starting from geometrical definition of this type of conoid, the cases $(i)$ and $(iii)$ will always have real plane sections.
\end{re}

\begin{thm}\label{Theorem_3_3}
If the solution of the system $(\ref{EQ_21})$ exists then it is not a non-degenerated conic.
\end{thm}

\noindent
{\bf Proof.} First we assume that $C \neq 0$. Then by substitution $z = -\mbox{\small $\displaystyle\frac{A}{C}$} \, x
- \mbox{\small $\displaystyle\frac{B}{C}$} \, y - \mbox{\small $\displaystyle\frac{D}{C}$}$ in~$(\ref{EQ_21}/i)$
we obtain $x-y$ projection of a possible intersection:
\begin{equation}
\label{EQ_28}
\begin{array}{l}
{\bigg (}
\mbox{\small $\displaystyle\frac{A^{2}\,a^{2}}{C^{2}}$}
+
\mbox{\small $\displaystyle\frac{A^{2}\,d^{2}}{C^{2}}$}
+
\mbox{\small $\displaystyle\frac {2\,d\,h\,A}{C}$}
{\bigg )}\,x^{2}\,y^{2}
+
{\bigg (}
\mbox{\small $\displaystyle\frac{2\,A\,B\,d^{2}}{C^{2}}$}
+
\mbox{\small $\displaystyle\frac{2\,d\,h\,B}{C}$}
+
\mbox{\small $\displaystyle\frac{2\,A\,B\,a^{2}}{C^{2}}$}
{\bigg )}\,x \, y^{3}
+                                                                                                 \\[1.65 ex]
{\bigg (}
\mbox{\small $\displaystyle\frac{B^{2}\,d^{2}}{C^{2}}$}
+
\mbox{\small $\displaystyle\frac{B^{2}\,a^{2}}{C^{2}}$}
{\bigg )}\,y^{4}
+
{\bigg (}
\mbox{\small $\displaystyle\frac{2\,A\,h\,a^{2}}{C}$}
+
\mbox{\small $\displaystyle\frac{2\,A\,D\,d^{2}}{C^{2}}$}
+
\mbox{\small $\displaystyle\frac{2\,A\,h\,d^{2}}{C} $}
+
\mbox{\small $\displaystyle\frac{2\,d\,h\,D}{C}$}
+
2\,d\,h^{2}
+                                                                                                 \\[1.65 ex]
\mbox{\small $\displaystyle\frac{2\,A\,D\,a^{2}}{C^{2}}$}
{\bigg )}\,x y^{2}
+
{\bigg (}
\mbox{\small $\displaystyle\frac{2\,B\,D\,d^{2}}{C^{2}}$}
+
\mbox{\small $\displaystyle\frac{2\,B\,D\,a^{2}}{C^{2}}$}
+
\mbox{\small $\displaystyle \frac {2\,B\,h\,a^{2}}{C}$}
+
\mbox{\small $\displaystyle\frac {2\,B\,h\,d^{2}}{C}$}
{\bigg )}\,y^{3}
+                                                                                                 \\[1.65 ex]
{\bigg (}
h^{2}\,b^{2}
-
\mbox{\small $\displaystyle\frac{A^{2}\,a^{2}\,b^{2}}{C^{2}}$}
{\bigg )}\,x^{2}
+
\mbox{\small $\displaystyle\frac{2\,A\,B\,a^{2}\,b^{2}}{C^{2}}$} \,x\,y
+
{\bigg (}
\mbox{\small $\displaystyle\frac{2\,D\,h\,a^{2}}{C}$} + h^{2}\,a^{2}
-
\mbox{\small $\displaystyle\frac{B^{2}\,a^{2}\,b^{2}}{C^{2}}$}
+                                                                                                 \\[1.65 ex]
\mbox{\small $\displaystyle\frac{D^{2}\,a^{2}}{C^{2}}$} + h^{2}\,d^{2}
+
\mbox{\small $\displaystyle\frac{2\,D\,h\,d^{2}}{C}$}
+
\mbox{\small $\displaystyle\frac{D^{2}\,d^{2}}{C^{2}}$}
{\bigg )}\,y^{2}
-
{\bigg (}
\mbox{\small $\displaystyle\frac{2\,A\,D\,a^{2}\,b^{2}}{C^{2}}$}
+
\mbox{\small $\displaystyle\frac{2\,A\,h\,a^{2}\,b^{2}}{C}$}
{\bigg )}\,x
-                                                                                                 \\[1.65 ex]
{\bigg (}
\mbox{\small $\displaystyle\frac{2\,B\,h\,a^{2}\,b^{2}}{C}$}
+
\mbox{\small $\displaystyle\frac{2\,B\,D\,a^{2}\,b^{2}}{C^{2}}$}
{\bigg )}\,y
-
h^{2}\,a^{2}\,b^{2}
-
\mbox{\small $\displaystyle\frac{2\,D\,h\,a^{2}\,b^{2}}{C}$}
-
\mbox{\small $\displaystyle\frac{D^{2}\,a^{2}\,b^{2}}{C^{2}}$}
=
0.
\end{array}
\end{equation}
The previous equation is algebraic equation of the second order iff:
$$
\begin{array}{l}
\mbox{\small $\displaystyle\frac{A^{2}\,a^{2}}{C^{2}}$}
+
\mbox{\small $\displaystyle\frac{A^{2}\,d^{2}}{C^{2}}$}
+
\mbox{\small $\displaystyle\frac {2\,d\,h\,A}{C}$}
= 0 \;\;\wedge                                                                                    \\[1.00 ex]
\mbox{\small $\displaystyle\frac{2\,A\,B\,d^{2}}{C^{2}}$}
+
\mbox{\small $\displaystyle\frac{2\,d\,h\,B}{C}$}
+
\mbox{\small $\displaystyle\frac{2\,A\,B\,a^{2}}{C^{2}}$}
= 0 \;\;\wedge                                                                                    \\[1.00 ex]
\mbox{\small $\displaystyle\frac{B^{2}\,d^{2}}{C^{2}}$}
+
\mbox{\small $\displaystyle\frac{B^{2}\,a^{2}}{C^{2}}$}
= 0 \;\;\wedge                                                                                    \\[1.00 ex]
\mbox{\small $\displaystyle\frac{2\,A\,h\,a^{2}}{C}$}
+
\mbox{\small $\displaystyle\frac{2\,A\,D\,d^{2}}{C^{2}}$}
+
\mbox{\small $\displaystyle\frac{2\,A\,h\,d^{2}}{C} $}
+
\mbox{\small $\displaystyle\frac{2\,d\,h\,D}{C}$}
+
2\,d\,h^{2}
+
\mbox{\small $\displaystyle\frac{2\,A\,D\,a^{2}}{C^{2}}$}
= 0 \;\;\wedge                                                                                    \\[1.00 ex]
\mbox{\small $\displaystyle\frac{2\,B\,D\,d^{2}}{C^{2}}$}
+
\mbox{\small $\displaystyle\frac{2\,B\,D\,a^{2}}{C^{2}}$}
+
\mbox{\small $\displaystyle \frac {2\,B\,h\,a^{2}}{C}$}
+
\mbox{\small $\displaystyle\frac {2\,B\,h\,d^{2}}{C}$}
= 0\,.
\end{array}
$$
Solving previous system, via Maple, we obtain two subcases:
$$
\mbox{\boldmath $1^{0}\!.$}\;
A\!=\!0, \, B\!=\!0, \, C\!=\!p, \, D\!=\!-ph
\quad \vee \quad
\mbox{\boldmath $2^{0}\!.$}\;
A\!=\!q, \, B\!=\!0, \, C\!=\!-\mbox{\small $\displaystyle\frac{a^2+d^2}{2dh}\,q$},
\, D\!=\!\mbox{\small $\displaystyle\frac{a^2+d^2}{2d}\,q$}
$$
for $p, q \in \mbox{\bf R} \backslash \{0\}$. Therefore:

\break

\noindent
\mbox{\boldmath $1^{0}\!.\;$}
\begin{minipage}[t]{121.0 mm}The system (\ref{EQ_21}):
$$
\left\{
\begin{array}{c}
(a^2 y^2 + d^2 y^2 - a^2 b^2)(z-h)^2 - 2dhxy^2(z - h) + b^2 h^2 x^2 = 0, \\[1.00 ex]
z - h = 0;
\end{array}
\right\}
$$
has the reduced Groebner basis (related to lexicographic order) in the following form:
$$
g_{11}:=x^2
\quad \wedge \quad
g_{12}:=z-h.
$$
Hence, the system (\ref{EQ_21}) has the intersection by directrix $d_{2}$ and this intersection is not non-degenerated conic.
\end{minipage}

\medskip
\noindent
\mbox{\boldmath $2^{0}\!.\,$}
\begin{minipage}[t]{121.0 mm}The system (\ref{EQ_21}):
$$
\left\{
\begin{array}{c}
(a^2 y^2 + d^2 y^2 - a^2 b^2)(z-h)^2 - 2dhxy^2(z - h) + b^2 h^2 x^2 = 0, \\[1.00 ex]
2dhx - (a^2+d^2)z + h(a^2+d^2) = 0;
\end{array}
\right\}
$$
has the reduced Groebner basis (related to lexicographic order) in the following form:
$$
g_{21}:=z^2 - 2zh + h^2
\quad \wedge \quad
g_{22}:=2dhx - (a^2+d^2)z + h(a^2+d^2).
$$
Hence, the system (\ref{EQ_21}) has intersection by directrix $d_{2}$ and this intersection is not non-degenerated conic.
\end{minipage}

\medskip
\noindent
Let $C=0$. We can consider two cases: $B = 0$ or $B \neq 0$. If $B = 0$ subcase $A = 0$ does not lead to a planar intersection. If $B = 0$
subcase $A \neq 0$ is considered in the Lemma 3.1. {\big (}case $(i)${\big )}. Finally, let us observe case $B \neq 0$.
By substitution of $y = -\mbox{\small $\displaystyle\frac{A}{B}$} x - \mbox{\small $\displaystyle\frac{D}{B}$}$
in~$(\ref{EQ_21}/i)$ we obtain $x-z$ projection of a possible intersection:
\begin{equation}
\label{EQ_29}
\begin{array}{l}
-
\mbox{\small $\displaystyle\frac{2\,d\,h\,A^{2}}{B^{2}}\;x^{3}$}\,z
+
{\bigg (}
\mbox{\small $\displaystyle\frac{a^{2}\,A^{2}}{B^{2}}$}
+
\mbox{\small $\displaystyle\frac{d^{2}\,A^{2}}{B^{2}}$}
{\bigg )}\;x^{2}\,z^{2}
+
\mbox{\small $\displaystyle\frac{2\,d\,h^{2}\,A^{2}}{B^{2}}\;x^{3}$}
-
{\bigg (}
\mbox{\small $\displaystyle\frac{4\,d\,h\,A\,D}{B^{2}}$}
+                                                                                                 \\[1.65 ex]
\mbox{\small $\displaystyle\frac{2\,h\,a^{2}\,A^{2}}{B^{2}}$}
+
\mbox{\small $\displaystyle\frac{2\,h\,d^{2}\,A^{2}}{B^{2}}$}
{\bigg )}\,x^{2}\,z
+
{\bigg (}
\mbox{\small $\displaystyle\frac{2\,a^{2}\,A\,D}{B^{2}}$}
+
\mbox{\small $\displaystyle\frac{2\,d^{2}\,A\,D}{B^{2}}$}
{\bigg )}\,x\,z^{2}
+
{\bigg (}
\mbox{\small $\displaystyle\frac{4\,d\,h^{2}\,A\,D}{B^{2}}$}
+                                                                                                 \\[1.65 ex]
\mbox{\small $\displaystyle\frac{h^{2}\,a^{2}\,A^{2}}{B^{2}}$}
+
h^{2}\,b^{2}
+
\mbox{\small $\displaystyle\frac{h^{2}\,d^{2}\,A^{2}}{B^{2}} $}
{\bigg )}\,x^{2}
-
{\bigg (}
\mbox{\small $\displaystyle\frac{2\,d\,h\,D}{B^{2}}$}
+
\mbox{\small $\displaystyle\frac{4\,h\,d^{2}\,A\,D}{B^{2}}$}
+
\mbox{\small $\displaystyle\frac{4\,h\,a^{2}\,A\,D}{B^{2}}$}
{\bigg )}\,x\,z
+                                                                                                 \\[1.65 ex]
{\bigg (}
\mbox{\small $\displaystyle\frac{a^{2}\,D}{B^{2}}$}
-
a^{2}\,b^{2}
+
\mbox{\small $\displaystyle\frac{d^{2}\,D}{B^{2}} $}
{\bigg )}\,z^{2}
+
{\bigg (}
\mbox{\small $\displaystyle \frac {2\,h^{2}\,d^{2}\,A\,D}{B^{2}}$}
+
\mbox{\small $\displaystyle \frac {2\,d\,h^{2}\,D}{B^{2}}$}
+
\mbox{\small $\displaystyle \frac {2\,h^{2}\,a^{2}\,A\,D}{B^{2}}$}
{\bigg )} \, x
+                                                                                                 \\[1.65 ex]
{\bigg (}
2\,h\,a^{2}\,b^{2}
-
\mbox{\small $\displaystyle \frac {2\,h\,a^{2}\,D}{B^{2}}$}
-
\mbox{\small $\displaystyle \frac {2\,h\,d^{2}\,D}{B^{2}} $}
{\bigg )} \, z
+
{\bigg (}
\mbox{\small $\displaystyle\frac{h^{2}\,d^{2}\,D}{B^{2}}$}
+
\mbox{\small $\displaystyle\frac{h^{2}\,a^{2}\,D}{B^{2}}$}
-
h^{2}\,a^{2}\,b^{2}
{\bigg )}
=
0.
\end{array}
\end{equation}
Therefore $A = 0$ and we obtain the case $(ii)$ from the Lemma 3.1.$\;$\stop

\begin{re}
\label{Remark_3_4}
Let us emphasize that we may consider possible factorizations of polynomials $(\ref{EQ_28})$ and $(\ref{EQ_29})$
in purpose of determining a type of section by non-degenerated conic as a part of plane intersection of this egg
curve based conoid.
\end{re}

\break

\section{\mbox{\normalsize \bf 4. CONCLUSIONS}}

In the paper an application of Groebner bases is considered, in order to define planarity of two surfaces intersection.
An egg curve based conoid is analyzed in particular, from the aspect of Projective Geometry. Its plane sections are analyzed
from the aspect of Analytic Geometry. Using the example of the egg curve based conoid, a method is developed, based on the technique
of Groebner bases, which enables investigation of existence of the assigned type of planar section. Based on the presented method,
in some cases, it is possible to determine even for the other algebraic surfaces if they would imply a planar section that consists
of non-degenerated conics.

\bigskip

\bigskip

\noindent
${}^{\ast}$Faculty of Electrical Engineering, Faculty of Civil Engineering; Belgrade; Serbia

\smallskip \noindent
$\;\,$E-mails: malesevic@etf.bg.ac.yu, marijao@grf.bg.ac.yu

\end{document}